\documentclass[12pt]{article} 
\usepackage{verbatim,amsmath,amssymb}
\usepackage{epsfig}
\usepackage{graphicx}
\usepackage{index}

\usepackage{amscd}
\input math111.def
\input buch.def
\numberwithin{equation}{section}

\def\span{\mbox{\rm span }}

\newif\ifmarglab
\makeatletter
\def
\label#1{\@bsphack\ifmarglab\marginpar{LAB:#1}\fi\if@filesw {\let\thepage\relax
   \def\protect{\noexpand\noexpand\noexpand}%
   \edef\@tempa{\write\@auxout{\string
      \newlabel{#1}{{\@currentlabel}{\thepage}}}}%
   \expandafter}\@tempa
   \if@nobreak \ifvmode\nobreak\fi\fi\fi\@esphack}
\makeatother

\begin{document}
\title{Two Remarks on Primary Spaces}
\author{Paul F. X . M\"uller
\footnote{Partly supported by FWF -project P 20166-N18}}
\date{October $18^{{\rm th}},$ 2009.}
\maketitle
\begin{abstract} 
We prove that for any operator $T$ on 
$ \ell^\infty( H^1 (\bT)) ,$ the identity factores 
through $ T $ or $ \Id - T .$ 

We  re-prove analogous  results of  H.M. Wark 
 for the spaces  $ \ell^\infty( H^p(\bT) ), $
 $1<p <\infty .$
In the present paper direct combinatorics of colored dyadic intervals 
replaces the dependence on Szemeredi's theorem in \cite{wa}.\\

\end{abstract}

\tableofcontents
\newpage
\section{Introduction}
Consider  $ \ell^\infty( H^p(\bT) ) ,$ consisting  of bounded 
$H^p (\bT)-$valued sequences,
and let $T$ be a continuous  linear operator acting on that space. 
The identity on  $ \ell^\infty( H^p (\bT)) ,$ factores  then 
through $ T $ or $ \Id - T .$ 
This assertion holds true for the range $ 1 < p \le \infty $ and also  
when $ H^p ( \bT ) $ is replaced by the dual space ${H^1 ( \bT )} ^* .$ 
See J. Bourgain \cite{b},  H.M. Wark \cite{wa}  and \cite{pfxm1,pfxm2}.
As a consequence  for each  $ 1 < p \le \infty $
the space $ \ell^\infty( H^p  ( \bT )) ,$ is  primary. 
 When  $p = \infty$ this class of spaces arises 
when the Pelczynski decomposition method is applied
to $H^\infty$ or its complemented subspaces. 
 %
%

The conceptual framework for proving these results was introduced by 
J . Bourgain~\cite{b}.  It consists of two basic steps and aims at 
replacing the factorization problem on $ \ell^\infty( H^p  ( \bT )) $
by  its localized, finite 
dimensional counterpart. The first step is to show that general operators on 
$  \ell^\infty( H^p ) ,$  can be replaced by  diagonal operators 
 on 
$$ \left( \sum_{n=1}^\infty  H^p_n  ( \bT ) \right) _\infty
. $$
The second (and often harder) part of the argument is to 
verify the following property
for operators acting on finite dimensional spaces. 
%
%
To each $n \in \bN $ and $ 1 < p \le \infty, $ there exists 
$ N = N( n , p ) $  so that for any linar operator on 
$H^p_N ( \bT) $,  the identity on $H^p_n ( \bT) $ 
can be well-factored  through $T $ or through $ \Id - T. $ 
That is,
\begin{equation}\label{factorization}
\begin{array}{ccc} H^p_n&\stackrel{\Id}{\longrightarrow}& H^p_n\\
E\downarrow&&\uparrow P\\
 H^p_N&\stackrel{H}{\longrightarrow}& H^p_N\end{array}
\end{equation} 
where $ H = T $ or $H = \Id - T $ 
and 
$ \|E \|\cdot \| P \| \le 16 . $

The methods yielding \eqref{factorization} fall roughly 
into two 
classes. The first one works is by {\it restriction:}
One  selects a 
large {\it subset} of a  basis in $H^p_N ( \bT) $
on which the operator acts as a multiplier. 
The second method uses averaging or {\it blocking :}
One    attempts to find  {\it block bases } 
 forming a set of approximate eigenvectors.
The three cases $ p = \infty ,$ $ 1 < p < \infty ,$  
and   ${H^1_n ( \bT )} ^* ,$ 
were treated by separate techniques as  
summarized below:
\begin{enumerate} 
\item The case ${p = \infty.}$ By harmonic analysis 
J. Bourgain \cite{b} proves that there exists $\lambda \in \bC $ and  an arithmetic progression $ \Lambda $
of length 
$ n$ in  $ \{ 1, \dots , N \} $ so that the restriction of $T $ to the 
exponentials $ e_j (t) = e^{i j t } $ satisfies
 $$ T (  e_j ) = \lambda   e_j  +\text{ small $L^\infty$ error} , 
\quad j \in   \Lambda . $$
\item The case  
 ${1 < p <  \infty.}$ 
As an application of Szemeredi's theorem \cite{sz},  H.M. Wark \cite{wa} 
obtains $\lambda \in \bC $ and 
an arithmetic progression $ \Lambda $
of length 
$ n$ in  $ \{ 1, \dots ,N \} $ so that the restriction of $T $ to the 
exponentials $ e_j (t) = e^{i j t } $ satisfies similarly 
 $ T (  e_j ) \sim  \lambda   e_j $  with  $j \in   \Lambda . $
Here the error is small in  $L^p.$

The $L^p$ boundedness of the Hilbert transform \cite{z} 
allows us to replace in \eqref{factorization} 
 the spaces of analytic polynomials  $H_n^p (\bT )$ 
respectively $H_n^p (\bT )$ 
finite dimensional Minkowski spaces $ \ell_m^p $
of the corresponding dimensions.
Now apply  J. Bourgain's and L. Tzafriri's results   on 
restricted invertibility. 
 To any  linear operator 
 $T : \ell_m^p \to  \ell_m^p $ 
 there exist  \cite{btz1}  random 
subsets $ \s$ of $ \{ 1 , \dots, m \} $ of cardinality
$ m /C $  so that the restriction of $T$ 
 to the corresponding 
$ \ell_m^p -$unit vectors   $\{x_j , j \in \s \}$ satisfies 
 $$ T (  x_j ) = \lambda_j   x_j  +\text{ small  $\ell_m^p $ error} , 
\quad\quad j \in   \s . $$ 
\item The case of   ${H^1_n ( \bT )} ^* .$
By  Bochkarev's theorem \cite[Section 4]{bo},  the spaces  ${H^1_n ( \bT )} ^* $
admit an  unconditional basis that is equivalent to the Haar basis 
in 
dyadic $BMO_n  .$  
Exploiting this connection \cite{pfxm2}, yields 
a  block basis $\{ b_I : I \in \cD _n \} $ in  ${H^1_N ( \bT )} ^* ,$ well-equivalent to the  
Bochkarev basis in  ${H^1_n ( \bT )} ^* $ on which the given operator $T :
{H^1_N ( \bT )} ^* \to {H^1_N ( \bT )} ^*$ acts as 
$$ T (  b_{J} ) = \lambda_J   b_J  +\text{ small $\BMO $ error} , 
\quad\quad J \in   \cD_n  . $$
\end{enumerate} 
In this note we adapt  the  blocking of the Haar basis
to the spaces 
$\ell^\infty(H^p ( \bT) ) ,$  $ 1\le p < \infty .$
First we settle the case of $\ell^\infty(H^1( \bT) ) ,$ left open in 
\cite{wa}.  After the reduction to diagonal operators, this  is obtained 
with  Bochkarev's theorem \cite[Section 4]{bo}, and
reduction to \cite[Theorem 2]{pfxm2}.
Second, we re-prove results of H. M. Wark \cite{wa}.
We find our block bases directly,  working with the Haar system; 
thus we provide an alternative to the method based on Szemeredi's theorem, 
and also    
to the restricted invertibility methods. 

\section{Notation and Preliminaries} 
\paragraph{The Haar System.}
We let $\cD $ denote the collection of (half-open) dyadic intervals
contained in the unit interval
$$
[(k-1)2^{-n} , k 2^{-n}[ , \quad 1\le k \le 2^{n} , \quad n \in \bN . $$
For $n \in \bN $ write  
$\cD_n = \{ I \in \cD : |I| \ge  2^{-n} \}. $  
Denote by $ \{ h_I :  I \in \cD \} $
the $L^\infty-$ normalized Haar system, where $h_I$ is supported on $I$ 
and 
$$
h_I = \begin{cases} 1 \quad &\text{ on the left half of }I ;\\ 
                    -1  \quad &\text{ on the right half of }I . 
      \end{cases}
$$ 
The Haar system is a Schauder basis in  $L^p,$ $(1 \le  p < \infty) $
and an unconditional Schauder basis when  $1 < p < \infty .$
For $ f \in  L^p$ define its dyadic square function as 
$$ S(f) = 
( \sum_{I \in \cD } \la f,\frac{h_I}{|I|}\ra ^2 1_I 
)^{1/2}.$$
Then,
$ c_p \| f\|_{L^p} \le \| S(f)\|_{L^p} \le C_p \| f\|_{L^p} , $
where $C_p \sim p^2/ (p-1)$ and $c_p  = C_p^{-1} ; $ see e.g., \cite[Chapter 1]{pfxm1}. 

To $ 1 \le p \le \infty $ and $ n \in \bN $ we denote by 
$ L^p_n $ the linear span of $ \{ h_I : I \in \cD_n\} , $
equipped with the norm in $ L^p . $
\paragraph{$H^p$ spaces.} Dyadic $H^1$ is defined by the relation 
$ f \in H^1$ if   $S(f) \in L^1$   and its norm is 
  $$  \| f\|_{H^1} = \| S(f) \|_{L^1} . $$
Let $ n \in \bN .$ 
We denote by $H^1_n(\d)$ the subspace of $H^1$ that is spanned by 
$\{ h_I : | I | \ge 2^{-n} \} $ endowed with the norm  in  $H^1 .$

Let $ \bT = \{ e^{i\theta} : \theta \in [-\pi,\pi] \}.$ 
For $1 \le p < \infty $ let $H^p(\bT)$ denote 
the norm closure in   $L^p(\bT)$ of 
$\span \{ e^{i n \theta } :  n \in \bN \} .$ 
The space $H^\infty(\bT)$ is the weak$-*$ closure in  $L^\infty(\bT)$
of $\span \{ e^{i n \theta } :  n \in \bN \} .$ 
As a consequence of the $L^p$ boundedness of the Hilbert transform \cite{z}
the spaces   $H^p(\bT)$ and  $L^p(\bT)$ are isomorphic. 
By  B. Maurey's theorem \cite{ma} the spaces   $H^1(\bT)$ and dyadic $H^1$ are 
isomorphic. 

For $ 1 \le p \le \infty$  we denote by   $H^p_{n }(\bT)$ 
the subspace of  $H^p(\bT)$ that is spanned by the polynomials
$ \{ e^{i\theta} , \dots , e^{i (2^{n+1} - 1)\theta} \} .$
The norm  on $H^p_{n }(\bT)$ is the one induced by  $H^p(\bT) .$
\paragraph{The Theorem of Bochkarev.} 
The trigonometric system does not form a Schauder basis in 
 $H^1_{n }(\bT)$ with basis constant independent of $n ,$ see \cite{z}.
Nevertheless the theorem of S.V. Bochkarev \cite{bo}
asserts that the spaces $H^1_{n }(\bT)$ admit an unconditional basis 
well equivalent to the Haar basis in 
  $H^1_n(\d).$ That is, there exists a system $\{ s_I \in  H^1_{n }(\bT)
: | I | \ge 2^{-n} \} $ 
of   functions 
so that the linear extension of the operator
$$ J : H^1_{n }(\bT) \to H^1_n(\d), s_I \to h_I $$
is bijective and defines an isomorphism with norm 
(and norm of the inverse) bounded independent of $n . $ 
Thus 
$$\| J : H^1_{n }(\bT) \to H^1_n(\d)\|\cdot\| J^{-1} : H^1_n(\d) \to  H^1_{n } (\bT)\| \le C .  $$
\paragraph{Diagonal operators.}
Let $Z_p = \left( \sum_{n=1}^\infty  H^p_n  ( \bT ) \right) _\infty. $
We say that $D : Z_p \to Z_p$  is a diagonal operator if
there exists a sequence of opertors 
$T_n :  H^p_n  ( \bT ) \to  H^p_n  ( \bT ) $ 
 so that 
$$ D( x_1, \dots , x_n, \dots) =  ( T_1(x_1) , \dots , T_n(x_n) , \dots) $$
\paragraph{Dyadic Trees.}
 Let $m \in \bN . $
We say that  $\{E_I: I \in \cD_m\}$  is  a dyadic tree of sets
if, the following conditions hold for every  $I \in \cD_{m-1} . $
If $I_1$  is the left half of  $I $ and 
 $I_2= I \sm I_1$   then 
\begin{equation}  E_{I_1} \cup   E_{I_2} \sbe   E_{I}  \quad\quad\text{and}   \quad\quad  E_{I_1} \cap   E_{I_2}
= \es . \end{equation} 

\paragraph{Lower $\ell^2$ estimates.}
 \cite[Chapter  5]{pfxm1}
 Let $ y _i \in H^1 $  , $ a _ i \in \bR $ 
and suppose that the sequence $ \{ y_i \} $ is disjointly supported over the 
Haar system. Then    
\begin{equation}\label{lowerl2}
 \left (\sum a_i^2 \|y_i \|^2 _{H^1}  \right)^{1/2} \le 
4  \|  \sum a_i^2 y_i \| _{H^1}  .
\end{equation}
\paragraph{The Carleson Constant.} Let $m \in \bN  $
and   $\cE = \{E_I: I \in \cD_m\}$    a dyadic tree.
 The Carleson constant  of 
$\cE$  is defined as  
\[ \y \cE \yy = \sup_{E_I\in \cE} \frac{1}{|E_I|} \sum_{E_J\in  \cE,\, E_J \sbe E_I } |E_J|. \] 
\paragraph{Projections and Large Carleson Constants.}
See \cite[Chapter 1, 3]{pfxm1}. Fix $N >> m , $ and let 
$ \cE _ I \sbe \cD _N,$ for 
$I \in \cD_m $
be  pairwise disjoint collections consisting of  pairwise disjoint
dyadic intervals,
so that 
$$ E_I = \bigcup_{J \in \cE_I } J \quad\text{and}\quad
b_I = \sum_{J \in \cE_I} h_J 
$$
satisfy the following conditions:
\begin{enumerate} 
\item The collection $ \{  E_I : I \in \cD_m \} $ is a dyadic tree satisfying
\begin{equation}\label{22okt1} 
\frac{|I|}{2}  \le | E_I|/|E_{[0,1[}| \le |I| .
\end{equation}
\item For each  $I \in \cD_{m-1} ,$ 
\begin{equation}\label{22okt2} 
E_{I_1} \sbe \{ t : b_I (t) = + 1 \} \quad \quad\text{and} \quad\quad
E_{I_2} \sbe \{ t : b_I (t) = - 1 \},\end{equation}
where  $I_1$ is   the left half of $I$ 
and $I_2= I \sm I_1.$
\end{enumerate}
Then \cite[Chapter 1]{pfxm1} the  orthogonal projection
\begin{equation}\label{22okt3} 
Q(f) = \sum_{J \in \cD_m } \la f , \frac{b_J}{\|b_J\|_2}
\ra  \frac{b_J}{\|b_J\|_2} , \quad \quad b_J = \sum_{I \in \cB_J} h_I ,
\end{equation}
is a bounded operator on $ L^ p $ and on $H^1$  with norm 
$ \|Q\|_p \le C_p ,$ respectively $   \|Q\|_{H^1} \le C .$

Let next $
\cL \sbe \{ E_I : I \in \cD_m    \} .
$
and
$$ Y_p = \left ( \span \{ b_I : I \in \cL \} , \| \,\,\| _p \right) . $$
By the condensation lemma in combination with the 
 Gamlen-Gaudet theorem  the following implication holds true,
see \cite[Chapter 3]{pfxm1}. 
If $n \in \bN $ is such that $ \y  \cL \yy \ge n 2^n ,$ then 
there exist linear operators 
$ E : L^p_n \to Y_p $
and
$R :  Y_p \to L_n^p $ 
so that 
\begin{equation}\label{22okt4}
\Id _{L^p_n} =    R  \Id _{Y_p}       I  
\quad\text{with}\quad  
 \|R\|_p  \cdot \| E \|_p   \le C_p . 
\end{equation}

\section{Factorization through Operators on 
$ \ell^\infty(H^1(\bT))$.}
In this section we prove that for any operator $T$ on 
$ \ell^\infty(H^1(\bT))$ the identity factores through
$T$ or $\Id - T . $ 
\begin{theor} \label{27okt1} Let $ T$ be a bounded operator on $ \ell^\infty( H^1(\bT) ) $
with $ \|T\| \le 1 . $ Then there exist an embedding $E$ and a projection $P$
so that  \begin{equation}\label{factorizationokt}
\begin{array}{ccc} \ell^\infty( H^1( \bT)) &\stackrel{\Id}{\longrightarrow}&  \ell^\infty( H^1( \bT))\\
E\downarrow&&\uparrow P\\
  \ell^\infty( H^1( \bT))&\stackrel{H}{\longrightarrow}&  \ell^\infty( H^1( \bT))\end{array} ,
\end{equation} 
where $ H = T $ or $ H = \Id - T $ and $ \|E\|\cdot \|P\| \le C . $
\end{theor}
Infer from  Maurey's theorem~\cite{ma}
that $H^1(\bT)$ can be replaced by dyadic $H^1 ,$ 
and from Wojtaszczyk's theorem \cite{woj10, woj20} that 
$\ell^\infty(H^1(\bT))$ is isomorphic to 
$ (\sum H^1_n(\d) )_\infty. $ 

First we give the  reduction of Theorem~\ref{27okt1} to diagonal operators on 
$$X= (\sum H^1_n(\d) )_\infty. $$  
\begin{theor} \label{27okt2} Let $ T$ be a bounded operator on $ \ell^\infty( H^1(\bT) ) $
with $ \|T\| \le 1 . $ Then there exist an embedding  
$ E : \ell^\infty( H^1( \bT)) \to X ,$  a projection $ P : X \to   \ell^\infty( H^1( \bT)) $ 
and a diagonal operator $D: X \to X$ satisfying $\|D\| \le C  $ and 
$$ D = PTE \quad\text{and (!)} \quad  \Id_X - D = P\left(\Id_{\ell^\infty( H^1(\bT) )}  - T\right)E, $$
where  $ \|E\|\cdot \|P\| \le C . $
\end{theor}
It is well known that the diagonal operator $D$ together with 
its embedding $E$ and projection $P$  can be constructed in straightforward manner 
from the assertions of the  
following proposition.
\begin{prop}\label{27okt3}
To $\e >0$ and $ n \in \bN $ there exists $ N = N( \e , n ) $ with the following property: To each  $n-$dimensional subspace $E \sbe H^1_N$ there exist pairwise disjount collections of dyadic intervals 
$ \{ \cB_J : J \in \cD_n \} $ so that 
\begin{enumerate}
\item
The sets $ B_J = \bigcup_{I \in \cB_J} I $ form a dyadic tree satisfying
$ |J|/2 \le |B_J |/|B_{[0,1[}  \le |J| .$
\item The orthogonal projection
$$ Q(f) = \sum_{J \in \cD_m } \la f , \frac{b_J}{\|b_J\|_2}
\ra  \frac{b_J}{\|b_J\|_2} , \quad \quad b_J = \sum_{I \in \cB_J} h_I ,$$
satisfies $\|Q\| \le 4 . $
\item For each $x \in E,  $ 
$$ \|Q(x)\|_{H^1} \le \e \|x\|_{H^1} . $$
\end{enumerate}
\end{prop}
\proof Let $\eta > 0 $ and $\{ x_1 , \dots, x_M \} $
 be an 
$ \eta -$net in $\{ x \in E : \|x \|_{H^1} = 1 \} .$
This can be done with $ M \le \exp ( 2n /\eta ) . $

Fix $ \tau = \tau( \e , n ) > 0 $ and put
$ \cL_i = \{ J \in \cD_N : |\la x , h_J \ra| \ge \tau \}. $
Apply \eqref{lowerl2} to obtain an upper  estimate for the cardinality of 
$  \cL_i . $
$$
\begin{aligned} 
\tau^{-1} \| x_i\|  &\ge  \|\sum _{J \in \cL _{i}}|J|^{-1}  h_J  \|_{H^1} \\
                    & \ge c | \cL_i | ^{1/2} , 
\end{aligned}
$$
for some $c>0 . $
Hence the union $ \cL = \bigcup_{i = 1 }^ M \cL_i $ is of cardinality 
$ \le Cn/\tau ^2. $
Put $ \cG_0  = \cD_N \sm   \cL . $ Its Carleson constant is bounded from below as 
$$ \y \cG_0 \yy \ge N - \ln ( c n \tau^{-2} ) . $$

In the next steps we select a subcollection $  \cG_M \sbe \cG_0 $ with the following two properties
\begin{enumerate} 
\item $\y \cG_M  \yy \ge \y \cG_0 \yy ( \tau ^{2M} c) $
\item For $ i \le M , $ 
$$ \int 
( \sum_{J \in \cG_M }\la  x_i , \frac{h_J}{|J|} \ra ^2 1_{J} 
)^{1/2} 
\le C \tau \|x_i \| _{H^1} .$$
\end{enumerate} 

The construction of $G_M$ is done inductively over $M$ steps.
In the first step we fix $x_1$ and let 
$ \{ J_k : 1 \le k \le N_0 \}$ with $ N_0 \le 2^{N - 1} $ 
be an enumeration of the intervals in $ \cG_0 . $
Define inductively a stopping time sequence $ n_0  , \dots ,n_L \le N_0 $  as follows.
 Put $n_o = 0  $ and 
$$ n_1 = \min \{ n \le N_0 : \int
( \sum_{k =  n_0 }^{n} \la  x_i , \frac{h_{J_k}}{|J_k|} \ra ^2 1_{J_k} 
)^{1/2} 
\le  \tau \|x_i \| _{H^1}\} .$$
Assume that $ n_0 , n_1 ,\dots , n_\ell $
 have been determined and that $ n_\ell < \infty . $     
Then define 
$$ n_{\ell +1} = \min \{ n \le N_0 : \int
( \sum_{k =  n_{\ell} +1  }^{n} \la  x_i , \frac{h_{J_k}}{|J_k|} \ra ^2 1_{J_k} 
)^{1/2} 
\le  \tau \|x_i \| _{H^1}\} .$$
Finally define
$$ L = \sup \{ \ell : n_\ell < \infty \} . $$  
We  use again \eqref{lowerl2}. Split $x_1$ along the Haar basis as
$ x_{1, \ell} = 
 \sum_{k =  n_{\ell} +1  }^{n_{\ell + 1 } } \la  x_1 , \frac{h_{J_k}}{|J_k|} \ra  h_{J_k} .
$
Note that $ x_{1, \ell}$ are disjointly supported over the Haar system 
so that the lower $\ell^2 $ estimate \eqref{lowerl2} gives
$$ 
\begin{aligned} 
\|x_1 \| _{H^1} &\ge \| \sum_{\ell = 1 } ^{L} x_{1, \ell} \|_{H^1} \\
                & \ge \tau \|x_1\|_{H^1} L^{1/2} . 
\end{aligned} 
$$   
Hence $ L \le C \tau^{-2} . $
Define next the partition of $ \cG_0 $ as 
$ \cG_0 = \bigcup_{\ell = 1 } ^L  \cG_{0, \ell} ,$ 
where $ \cG_{0, \ell} =  \bigcup_{k =  n_{\ell} +1  }^{n_{\ell + 1 } } \{J_k\} . $
Consequently there exists $ \ell_0 \le L $ so that
$  \y  \cG_{0, \ell_0} \yy \ge  \y \cG_0 \yy /L . $ Put 
$ \cG_1 =  \cG_{0, \ell_0} . $
Summing up the first step, we  obtained $ \cG_1\sbe  \cG_0 $ so that 
$$ \y \cG_1 \yy \ge  c\tau^2 \y \cG_0 \yy  , $$
and
$$ \int 
( \sum_{J \in \cG_1 }\la  x_1 , \frac{h_J}{|J|} \ra ^2 1_{J} 
)^{1/2} 
\le C \tau \|x_1 \| _{H^1} .$$
This completes the first step. In the second step 
we repeat the argument of the first step with 
$ x_ 1 $ replaced by $ x_2 $ and $\cG_0 $ replaced by   $\cG_1 .$
This implies the existence of $ \cG_2 \sbe \cG_1 $ so that 
$$ \y \cG_2 \yy \ge  c\tau^2 \y \cG_1 \yy  , $$
and
$$ \int 
( \sum_{J \in \cG_2 }\la  x_i , \frac{h_J}{|J|} \ra ^2 1_{J} 
)^{1/2} 
\le C \tau \|x_i \| _{H^1} \quad\text{for}\quad i \in \{ 1,2\} .$$
Iterating $ M -$times we obtain a decreasing chain of collections
$ \cG_M \sbe \dots \sbe \cG_0 , $ so that 
$$ \y \cG_M \yy \ge  c\tau^{2M} \y \cG_0 \yy  , $$
and 
$$ \int 
( \sum_{J \in \cG_M }\la  x_i , \frac{h_J}{|J|} \ra ^2 1_{J} 
)^{1/2} 
\le C \tau \|x_i \| _{H^1} \quad\text{for}\quad 1 \le i \le M  .$$
We specify  now $ \tau , \eta $ so that 
 $ ( \tau + \eta ) < \e 2^{-n} ( \log n )^{-1} . $
Consequently there exists $ N = N( \e , n )$ so that 
$$ \y \cG_M \yy \ge  c\tau^{2M} \cdot N \ge 4^n  . $$
Next we apply the condensation lemma and the Gamlen-Gaudet 
selection process \cite[Chapter 3]{pfxm1} to $  \cG_M . $

This implies that there exist $ \{ \cB_J \sbe \cG_M : J \in \cD_n \} $ 
so that  $ B_J = \bigcup_{I \in \cB_J} I $ and $b_J = \sum_{I \in \cB_J}h_I $
 satisfy the following conditions:
\begin{enumerate} 
\item  $\{ B_J  : J \in \cD_n \} $ is a dyadic tree and 
        $$ |J|/2 \le |  B_J |/|B_{[0,1[}| \le |J| . $$ 
\item $ B_{J_1} \sbe \{t:  h_I ( t) = +1 \}$ and $ B_{J_2} \sbe \{t:  h_I ( t) = -1 \},$
where $J_1 $ is the left half of $J$ and  $J_2= J \sm J_1 .  $
\end{enumerate} 
Next we claim that moreover
  $$|\la x , b_J \ra| \le ( \tau + \eta ) \|x\|_{H^1}\quad x \in E$$
%
To see this fix $x \in E $ with $\|x\|_{H^1}= 1$
Choose $x_i \in \{ x_1 , \dots, x_m \} $  so that $\|x- x_i \|_{H^1}< \eta ,$
and estimate using the defining properties of $ \cG_M . $
$$ 
\begin{aligned}
|\la x , b_J \ra| & =  |\la x - x_i  , b_J \ra| +| \la x_i  , b_J \ra|\\
                  & \le (\eta + \tau) .
\end{aligned} 
$$
It remains to verify that $ \|Q(x)\|_{H^1} \le \e 
$ for 
$ x \in E . $
Using the above coefficient estimate and triangle inequality we get,
$$ 
\begin{aligned} \|Q(x)\|_{H^1} & \le ( \tau + \eta) 
\| \sum_{J \in \cD_n} b_J /\|b_J\|_2^2 \|_{H^1} \\
& \le  ( \tau + \eta)  2^n \ln (n)  
.
\end{aligned}
$$
\endproof

We have now reduced the problem of factoring  the identity on 
$ \ell^\infty(H^1(\bT))$ to its finite dimensional counterpart.
The next step in the proof is 
factorization of the identity on $ H^1_{N }(\bT) . $ 
The following theorem settles a point left open in \cite{wa}.
\begin{theor}
For any $ n \in \bN $ there exists
$ N = N ( n ) $ such that for any linear operator
$ T :  H^1_{N }(\bT) \to H^1_{N } (\bT) $ 
the identity on $ H^1_{n } (\bT) $
factores through $ H = T $ or $ H = \Id_{ H^1_{N }(\bT) } - T, $ as 
$$\begin{array}{ccc}{ H^1_{n }(\bT)}&\stackrel{\Id}{\longrightarrow}&{ H^1_{n }(\bT)}\\
E\downarrow&&\uparrow P\\
{ H^1_{N} (\bT)}&\stackrel{H}{\longrightarrow}&{ H^1_{N }(\bT)}\end{array}$$
where $ E  : H^1_{n } (\bT) \to H^1_{n } (\bT) $ and 
$ P :  H^1_{n } (\bT) \to  H^1_{n } (\bT) $ are bounded linear operators saitsfying
$$ \|E \|\cdot \| P \| \le C . $$
\end{theor} 
\proof
Apply the theorem of S.V. Bochkarev \cite[Section 4]{bo}, 
and dualize  \cite[Theorem 2]{pfxm2}.
\endproof

\section{Factorization through Operators on 
$\ell^\infty (L^p),$ $1 < p < \infty.$}
In this section  we re-prove the theorem of H.M. Wark \cite{wa}  that
for  any bounded operator $T$ on $\ell^\infty (L^p)$ the identity factores 
through $T$ or $\Id - T . $ 
\begin{theor} \label{37okt1} Let $ T$ be a bounded operator on $ \ell^\infty( L^p(\bT) ) $
with $ \|T\| \le 1 . $ Then there exist an embedding $E$ and a projection $P$
so that  \begin{equation}\label{factorizationokt}
\begin{array}{ccc} \ell^\infty( L^p( \bT)) &\stackrel{\Id}{\longrightarrow}&  \ell^\infty( L^p( \bT))\\
E\downarrow&&\uparrow P\\
  \ell^\infty( L^p( \bT))&\stackrel{H}{\longrightarrow}&  \ell^\infty( L^p( \bT))\end{array} ,
\end{equation} 
where $ H = T $ or $ H = \Id - T $ and $ \|E\|\cdot \|P\| \le C_p . $
\end{theor}
It is well known that Theorem~\ref{37okt1} 
follows   from the next theorem.
\begin{theor}\label{Lp}
For any $ 1 < p < \infty $  and any $ n \in \bN $ there exists
$ N = N ( n, p ) $ so that for any linear operator
$ T : L^p_N \to L^p_N $ the identity on $ L^p_n $
factores through $ H = T $ or $ H =\Id_{ L^p_N} - T, $ as 
$$\begin{array}{ccc} L^p_n&\stackrel{\Id}{\longrightarrow}& L^p_n\\
E\downarrow&&\uparrow P\\
 L^p_N&\stackrel{H}{\longrightarrow}& L^p_N\end{array}$$
where $ E  : L^p_n \to L^p_N $ and $ P : L^p_N \to L^p_n $ are bounded linear operators satisfying
$ \|E \|\cdot \| P \| \le C_p . $
\end{theor} 
We obtain the finite dimensional factorization of Theorem~\ref{Lp} by exhibiting a 
block-basis of the Haar system 
on which $T$ acts as a multiplier. We use  combinatorics of colored dyadic intervals as in \cite{pfxm1}. 
(For an alternative derivation  see the appendix.) 
The desired block bases are constructed 
below. 
\begin{theor} \label{gg} For every $ 1 < p < \infty , $ and $
m \in \bN $ there exists $ N = N (m ,p ) $ so that the following 
holds true: 
To every linear operator $ H  : L^p_N \to L^p_N $  bounded by 
 $ \|H \| \le 1 $ 
there exist pairewise disjoint collections 
$$ \cE _ I \sbe \cD _N,\quad I \in \cD_m ,$$
consisting 
of pairwise disjoint intervals so that 
\begin{equation}\label{22okt11} 
E_I = \bigcup_{J \in \cE_I } J \quad\text{and}\quad
b_I = \sum_{J \in \cE_I} h_J ,
\end{equation}
satisfy the following conditions:
\begin{enumerate} 
\item The collection $ \{  E_I : I \in \cD_m \} $ is a dyadic tree satisfying
\begin{equation}\label{22okt12} 
 \frac{|I|}{2}  \le | E_I| \le |I| .
\end{equation}
\item For each  $I \in \cD_{m-1} ,$ 
\begin{equation}\label{22okt13} 
 E_{I_1} \sbe \{ t : b_I (t) = + 1 \} \quad \quad\text{and} \quad\quad
E_{I_2} \sbe \{ t : b_I (t) = - 1 \},
\end{equation}
where  $I_1$ is   the left half of $I$ 
and $I_2=I\sm I_1 .$
 \item For each $I \in \cD_m ,$ 
\begin{equation}\label{22okt14} 
\sum_{\{J \in \cD_m : J \ne I \}} |\la Hb_J, b_I \ra|  
\le |I|^{4} . 
\end{equation}
\end{enumerate} 
\end{theor} 
The following Lemma~\ref{thin} is the  main component in the proof of Theorem~\ref{gg}. 
Fix $ 1 < p < \infty , $  $I \in \cD $  and $ k \in \bN .$ 
Let   $ x \in L^p , $  $ y \in L^q $  and $1/p + 1/q = 1 $ so that 
\begin{equation}\label{10okt1}
 \|x\|_p\le |I| ^{1/p} \quad\text{and}\quad \|y \|_q \le  |I| ^{1/q}.
\end{equation} 
Define $ \cB$ to be  the collection of dyadic intervals for which the Haar coefficients of 
$ x $ resectively $ y $ are above the critical  threshold $ |J|/k ,$
thus 
$$ \cB = \left\{ J \in Q(I) : |\la x, h_J \ra| + | \la y , h_J \ra | > 
|J|/k \right\} . $$
\begin{lemma}\label{thin}
Let $ \ell \in \bN $ and put 
$$
A_p = ( k^2\ell^2 ) ( C_p + C_q ) + 1 .$$
There exists $ 1\le j \le A_p $ so that the collection 
$ \cB_j = \{ J \in \cB : |J| = 2^{-j} |I| \} $ 
satisfies
$$
\sum_{J \in \cB_j}  |J| \le \frac{|I|}{\ell} 
$$
\end{lemma}
\proof
Assume that the lemma is false. 
That is, for each $ j \in A_p $ we have 
\begin{equation} \label{16okt2}
\sum_{J \in \cB_j}  |J| > |I|/\ell.
\end{equation}
For $ J \in \cB_j $
there holds 
$ |\la x, h_J \ra| + | \la y , h_J \ra | > 
|J|/k . $ 
Hence by summing the lower estimate \eqref{16okt2} we obtain 
\begin{equation} \label{16okt3}
\sum _{j = 1} ^{ A_p}   \sum_{J \in \cB_j} 
|\la x, h_J \ra| + | \la y , h_J \ra | \ge A_p \frac{|I|}{k\ell}.
\end{equation}
On the other hand it follows from the unconditionality of the Haar system 
in $ L^p , 1 < p < \infty $ 
that for any choice of $ \a_J \in [-1, +1 ] $ 
\begin{equation} \label{16okt4}
\|\sum _{j = 1} ^{ A_p}   \sum_{J \in \cB_j} 
\a_J h_J \|_{L^p } \le C_p \sqrt{ A_p } |I|^{1/p} , \end{equation}
Invoking the assumptions that  $ \|x\|_p\le |I| ^{1/p}$
and $ \|y \|_q \le  |I| ^{1/q}$, the unconditionality of the Haar system
\eqref{16okt4} yields 
the following upper estimate for the Haar coefficients
competing with \eqref{16okt3}, 
\begin{equation} \label{16okt5}
\begin{aligned}
\sum _{j = 1} ^{ A_p}   \sum_{J \in \cB_j} 
|\la x, h_J \ra| + | \la y , h_J \ra | & \le  C_q \sqrt{ A_p } |I|^{ 1/p + 1/q}  + C_p  \sqrt{ A_p } |I|^{ 1/q + 1/p}\\
&\le  (C_q +C_p) \sqrt{ A_p }  |I| .
\end{aligned}
\end{equation}
Comparison of \eqref{16okt3} with \eqref{16okt5}
gives an upper estimate for $ A_p $ as follows,
$$ A_p \frac{|I|}{k\ell} \le  (C_q +C_p) \sqrt{ A_p }  |I| .$$
By cancellation and arithmetic we get 
$ A_p \le \ell^2 k^2 (C_q +C_p)^2 , $
contradicting the initial choice of $ A_p . $
\endproof

\paragraph{Proof of Theorem \ref{gg}.}
 Put $ b_1 = h_{[0,1]} , $ and $ \cE_1 = \{ [0,1[\} . $ At stage $ i $ 
we are given pairwise disjoint collections of dyadic intervals in $
\cD _N ,$ $ \cE_1, \dots , \cE_i $ and functions
$ b_1 , \dots , b_N  $ so that 
\begin{equation} \label{16okt6}
\sum _{ j =1}^{ i-1} |\la H b_j , b_i \ra| + | \la H^* b_j , b_i \ra| \le \|b\|_2^2 4^{-i} . \end{equation}
Fix $ \cE_i $ and $ J \in Q(\cE_i).$
Put 
$$g(J) = \sum _{
\{j : j \le i\}} |\la H b_j , h_J  \ra| + | \la H^* b_j , h_J \ra|
$$
and define the collections of good intervals as follows
$$
\cG = \left\{ J \in Q(\cE_i) : g(J)  \le |J| 4^{-i-1} \right\}. $$ 
By  Lemma~\ref{thin}  there exists 
             \begin{equation} \label{27okt11}
\mu \le 2^{(C_p + C_q)^2 (i+1)}, 
\end{equation}
so that the collection of dyadic intervals 
$$ 
\cC = \bigcup_{I \in \cE_i} \{ J \in \cG : J \sbe I , |J| = 2^{-\mu} |I| \} , $$
satisfies the lower estimate 
\begin{equation} \label{16okt8}
\sum_{J \in \cC}  |J| \ge ( 1 - 8^{-i} ) \sum_{J \in \cE_i}  |J| . 
\end{equation}
Next we define the collection of dyadic intervals $ \cE_{i+1} $ and the associated 
function $b_{i+1} $ using  to the well known
Gamlen Gaudet construction.  In case  $i$ is even we put
$$
\cE_{i+1} =
           \left\{ I \in \cC : I \sbe \{t : b_{\frac{i}{2}(t)} = -1\}\right\} ,
\quad b_{i+1} = \sum _{ J \in \cE_{i+1}} h_J .$$
Similarly if  $i $ is odd, we put 
  $$\cE_{i+1}    =     
\left\{ I \in \cC : I \sbe \{t : b_{\frac{i-1}{2}(t)} = +1\}\right\} .  
$$
Since $ \cE_{i+1} \sbe \cC \sbe \cG $   our construction so far, 
implies that 
$b_{i+1}$
satisfies the almost orthogonality relation 
\begin{equation} \label{16okt7}
\sum _{j = 1 }^{ i} |\la H b_j , b_{i+1} \ra| + | \la H^* b_j , b_{i+1} \ra| \le \|b_{i+1}\|_2^2 4^{-i-1} . 
\end{equation} 
In view of \eqref{16okt6} and \eqref{16okt7} we may now repeat the induction 
step. Taking into account the (cummulative effect of) the  lower estimate
\eqref{16okt8} and \eqref{27okt11} it is easy to see that we may choose $N = N(\e , n ) $ large enough
so that we may continue iterating  
until we reach 
 step $ 2^{m+1} - 1. $ 
The system
$ \{ b_1 , \dots ,b_{ 2^{m+1} - 1} \} $ satisfies then the full set of 
almost orthogonalty relations 
$$
\sum_{\{j : j \ne i\} } |\la Hb_j, b_i \ra|  
\le \|b_i\|_2^2 4^{-i} , \quad i \le 2^{m+1} - 1 .  $$
It remains to relabel the collections $\{\cE_{i} \} $ and functions
$\{b_{i}\}$ appropriately.
Since $ \mu \le m $ and $ 1 \le k \le 2^{\mu} $ define uniquely 
$$ j = 2^{\mu } +k - 1 \quad\text{and}\quad J = \left[ \frac{k-1}{2^\mu} , \frac{k}{2^\mu} \right[, $$
we obtain a canonical bijection between $ \cD_m $ and $ \{ 1 ,\dots , 2^{m+1} - 1\}. $
We relabel accordingly using this correspondence,
$$ \cE_{J} = \cE_j , \quad \quad b_J = b_j . $$
As a result of the Gamlen Gaudet construction the sets
$$ E_J = \bigcup _{ J \in \cE_{J} } J , \quad J \in \cD_m ,$$
form a dyadic tree, and 
the conditions \eqref{16okt8} translate into a measure estimate
$$ |J|/2 \le | E_J| \le |  J | ,\quad J \in \cD_m .$$
By \eqref{22okt1}---\eqref{22okt3},  the orthogonal projection
$$ Q(f) = \sum_{J \in \cD_m } \la f , \frac{b_J}{\|b_J\|_2}
\ra  \frac{b_J}{\|b_J\|_2} , \quad \quad b_J = \sum_{I \in \cB_J} h_I$$
is bounded operator on $ L^ p $ with norm $ \|Q\|_p \le C_p .$ 
\endproof

\paragraph{Proof of Theorem \ref{Lp}.}
Fix $ n \in \bN $ and put
$ m = C_p n2^n . $
Determine $ N = N(m, p ) $ so that the conclusion of 
Theorem~\ref{gg}  holds true. Accordingly we may select 
pairewise disjoint collections 
$ \{\cE _ I \sbe \cD _N : I \in \cD_m\} $
 consisting 
of pairwise disjoint intervals so that conditions 
\eqref{22okt11}--\eqref{22okt14} hold true. 
%

Next define the collections
$$
\cL = \left\{ E_I : | \la Tb_I , b_I \ra | \ge \|b_I\|^2_2/2 \right \}
$$
and  $\cR = \cL\sm \{ E_J :  J \in \cD_m \} .$
Clearly, then  
$ 
\cL \cup \cR = \{ E_J :  J \in \cD_m \} ,$ and hence,
one of these collections satisfies a lower estimate for its Carleson Constant, 
$$ \y  \cL \yy \ge m/2 \quad \text{or} \quad \y  \cR \yy \ge m/2 . $$
We continue under the assumption that the first alternative holds
and remark that otherwise we would replace $T$ by $\Id - T$ and continue
with the argument given below. 
Define, then, the space 
$$ Y_p = \left ( \span \{ b_I : I \in \cL \} , \| \,\,\| _p \right) . $$
Since $ \y  \cL \yy \ge n 2^n/2 $, \eqref{22okt1}--\eqref{22okt4}
imply that  there exist linear operators 
$$ I : L^p_n \to Y_p \quad\quad\text{and}\quad\quad  R :  Y_p \to L_n^p $$ 
so that 
$\Id _{L^p_n} =    R  \Id _{Y_p}       I  $
and 
$ \|R\|_p  \cdot \| I \|_p   \le C_p. $
On  $L^p_n$  define the operator
$$ P(f) = \sum_{J \in \cD_m } \la f , b_J
\ra  b_J \la T b_J , b_J \ra ^{-1}  , $$
By choice of $ \cL $ we get with \eqref{22okt1}--\eqref{22okt3} that 
$$\| P(f) \|_p \le 2C_p \| Q(f) \|_p \le 2C_p \| f \|_p . $$
We will verify next  that for $ g \in Y_p ,$ 
\begin{equation}\label{22okt21} 
\|P T g - g \|_{p} \le \frac12 \|g \|_p . 
\end{equation}
To this end fix  $I \in \cL $  and expand $ PT(b_I)$ as 
$$ PT(b_I) = 
b_I +  
 \sum_{ \{J \in \cL : J\ne I \}}  \la Tb_J , b_I \ra 
b_J  \la T b_J , b_J \ra ^{-1}.
$$
Next let  $  g \in Y_p ,$  hence 
$ g = \sum_{I \in \cL }c_I b_I .$
Put 
$$ \b(J) = 
\sum_{ \{I \in \cL : J\ne I \}}  c_I \la Tb_J , b_I \ra 
$$
and define the error term 
$$  e =\sum_{\{ J \in \cL \}}  \b(J)
b_J  \la T b_J , b_J \ra ^{-1} . $$
Then, by arithmetic we get 
$$ PTg =  g + e .$$
Since 
$
\b(J)
\le C_p^{-1} |J|^4 \|g\|_p , $
and 
$ |\la T b_J , b_J \ra| \ge \|b_J\|_2^2/2 ,$
whenever 
$J \in \cL  , $
we obtain with triangle inequality, the following  estimate for the off-diagonal term,
$$ \| e \| _p \le \|g\|_p /2 , $$
which gives \eqref{22okt21}.
\endproof

\paragraph{Appendix: Reviewing  restricted invertibility.} A proof of 
Theorem~\ref{Lp}  follows directly from the random methods 
of J. Bourgain and L. Tzafriri \cite{btz1}
on restricted invertibility of matrices.
\begin{theor}
For any $ 1 < p < \infty $  and any $ n \in \bN $ there exists
$ N = N ( n, p ) $ 
so that for any linear operator
$ T : \ell^p_N \to \ell^p_N $ the identity on $ \ell^p_n $
factores through $ H = T $ or $ H =\Id_{ \ell^p_N} - T, $ as 
$$ \Id_{ \ell^p_n} = P H E $$
where $ E  : \ell^p_n \to \ell^p_N $ and $ P : \ell^p_N \to \ell^p_n $ are bounded linear operators satisfying
$$ \|E \|\cdot \| P \| \le C . $$
\end{theor} 
\proof The operator $ T $ has matrix representation with respect to 
$ \{ e_j : j \le n \} $, the unit vector basis of $ \ell^p_N .$
Define  $ \cL = \{ i \le n : |( Te_i, e_i )| \ge 1/2 \} $ and 
$ \cR = \{ 1 , \dots , n \} \sm \cL . $ Thus  the cardinality of
$ \cL $ or $ \cR $ is at least $ N/2. $ Assume without loss of generality
that $| \cL | \ge N/2 . $ (Otherwise replace 
$T$ by $\Id - T $ and proceed as below.) 
Put   $ N_1 = |\{ 1 , \dots , N\} \cap \cL |.$ 
Clearly, there exist now linear operators
$E_1 : \ell^p_N \to \ell^p_{N_1} $ and  $P_1  : \ell^p_{N_1} \to \ell^p_N ,$
and an operator $ S : \ell^p_{N_1} \to \ell^p_{N_1} $
satisfying the following  conditions:
First, $S$  factores through $T$ as 
$$ S = P_1 T E_1 , \quad\text{and} \quad \|E_1 \| \cdot \| P_1 \| \le C $$
and, second,
$$ ( Se_i , e_i ) = 1 \quad \text{for} \quad 1 \le i \le N_1. $$
Now we apply the theorem of J. Bourgain and L. Tzafriri 
\cite{btz1}  on restricted invertibility to the operator 
$ S  : \ell^p_{N_1} \to \ell^p_{N_1} .$
To each given $ \e > 0 $ 
it yields the existence of a constant $ C = C( \e , p ) $ and subset 
$ \s \sbe \{ 1 , \dots , N_1\}$ so that $ | \s | \ge N_1 /C $ and 
$ R_\s S R_s $ is invertible (on its range) and satisfies 
$$ \|  R_\s S R_\s \|_p \le 1 + \e . $$
(Here $R_\s$ denotes the operator of restriction  to basis vectors 
$ \{ e_j : j \in \s \} .$ ) 
Now put $ N_2 = |\{ 1 , \dots , N_1\} \cap \s | .$ 
By restricted invertibility there are now operators 
$E_2 : \ell^p_{N_2} \to \ell^p_{N_1} $ and  
$P_1  : \ell^p_{N_1} \to \ell^p_{N_2} ,$ so that the identity 
on $\ell^p_{N_2}$ factores through $S$ as 
$$ \Id_{\ell^p_{N_2}} = P_2 S E_2  \quad\text{with} \quad \|E_1 \| \cdot \| P_1 \| \le 2 .$$
In summary, we obtained the factorization through $T ,$
$$ \Id_{\ell^p_{N_2}} = P_2 P_1T E_1E_2  
\quad\text{with } N_2 \ge N/2C . $$

\bibliographystyle{abbrv}
\bibliography{primary}

\paragraph{Address:} Institut f\"ur Analysis \\
J. Kepler Universit\"at\\
A-4040 Linz\\
Austria.\\
pfxm@bayou.uni-linz.ac.at
\end{document}